\documentclass{amsart}
\theoremstyle{definition}

\theoremstyle{remark}

\numberwithin{equation}{section}

\begin{document}

\title{Diameter of the commuting graph of $\mathbb{R}^{n\times n}$}

\author{Yaroslav Shitov}
\address{National Research University Higher School of Economics, 20 Myasnitskaya Ulitsa, Moscow 101000, Russia}
\email{yaroslav-shitov@yandex.ru}


\begin{abstract}
The vertices of \textit{commuting graph} of $\mathbb{R}^{n\times n}$ are non-scalar matrices;
the edges are defined as pairs $(u,v)$ satisfying $uv=vu$. One can see~\cite{AMRP} that,
for $n\geq 3$, the diameter of this graph is at least four; we give a short proof that it is exactly
four. A longer proof of this result has been known~\cite{Grau, Miguel}.
\end{abstract}

\maketitle
If the eigenvalues of $a$ are simple, then $\mathbb{R}^n$ is a direct sum of invariant subspaces
each of which has dimension either $1$ or $2$; in this case, $a$ commutes with a nontrivial idempotent.
The distance between a pair of idempotents is at most $2$ by Theorem~11 of~\cite{AMRP},
so that, for any matrices $a_0$ and $a_4$, there are sequences $a_0(t),\ldots,a_4(t)$ of non-scalar matrices such that
$a_0(t)$ and $a_4(t)$ converge to $a_0$ and $a_4$, respectively, and $a_{j-1}(t)$ commutes with $a_{j}(t)$.
Multiplying $a_i(t)$ by a constant and adding a scalar matrix,
we can get a singular matrix $\alpha_i(t)$ with unit norm.
Since the unit sphere is compact, every $\alpha_i(t)$ is dense in some $\alpha_i$, which is non-scalar
and real up to adding a scalar matrix. Passing to the limit, we see that $\alpha_{j-1}$ and $\alpha_{j}$ commute.

\end{document}